\documentclass[reqno,11pt]{amsart} 
\usepackage{amsmath}
\usepackage{amssymb}
\usepackage{amsthm}
\usepackage{verbatim}
\usepackage{xcolor}
\usepackage{graphics}

\newcommand\NoBlackBoxes{\global\overfullrule0pt}

\NoBlackBoxes
\theoremstyle{plain}

\begin{document}

\title{Strictly subgaussian probability distributions
}

\author{S. G. Bobkov$^{1,4}$}
\thanks{1) 
School of Mathematics, University of Minnesota, Minneapolis, MN, USA,
bobkov@math.umn.edu. 
}

\author{G. P. Chistyakov$^{2,4}$
}

\author{F. G\"otze$^{2,4}$}
\thanks{2) Faculty of Mathematics,
Bielefeld University, Germany,
goetze@math-uni.bielefeld.de.
}
\thanks{4) Research supported by the NSF grant DMS-2154001 and 
the GRF – SFB 1283/2 2021 – 317210226 }
\dedicatory{Dedicated to the memory of Gennadiy P. Chistyakov\;  *\,May 1, 1945 \; $\dagger$\,December 30, 2022.} 


\begin{abstract}
We explore the class of probability distributions on the real line whose Laplace 
transform admits a strong upper bound of subgaussian type. Using Hadamard's 
factorization theorem, we extend the class $\mathfrak L$ of Newman and propose new sufficient 
conditions for this property in terms of location of zeros of the associated characteristic 
functions in the complex plane. The second part of this note deals with Laplace 
transforms of strictly subgaussian distributions with periodic components. 
This subclass contains interesting examples, for which the central limit theorem 
with respect to the R\'enyi entropy divergence of infinite order holds.
\end{abstract}

\maketitle
\markboth{S. G. Bobkov, G. P. Chistyakov and F. G\"otze}{Strictly subgaussian distributions}

\def\theequation{\thesection.\arabic{equation}}
\def\E{{\mathbb E}}
\def\R{{\mathbb R}}
\def\C{{\mathbb C}}
\def\P{{\mathbb P}}
\def\Z{{\mathbb Z}}
\def\S{{\mathbb S}}
\def\I{{\mathbb I}}
\def\T{{\mathbb T}}

\def\s{{\mathbb s}}

\def\G{\Gamma}

\def\Ent{{\rm Ent}}
\def\var{{\rm Var}}
\def\Var{{\rm Var}}
\def\V{{\rm V}}

\def\H{{\rm H}}
\def\Im{{\rm Im}}
\def\Tr{{\rm Tr}}
\def\s{{\mathfrak s}}

\def\k{{\kappa}}
\def\M{{\cal M}}
\def\Var{{\rm Var}}
\def\Ent{{\rm Ent}}
\def\O{{\rm Osc}_\mu}

\def\ep{\varepsilon}
\def\phi{\varphi}
\def\vp{\varphi}
\def\F{{\cal F}}

\def\be{\begin{equation}}
\def\en{\end{equation}}
\def\bee{\begin{eqnarray*}}
\def\ene{\end{eqnarray*}}

\thispagestyle{empty}

\vskip-10mm
\section{{\bf Introduction}}
\setcounter{equation}{0}

\vskip2mm
\noindent
Following Kahane \cite{Kah}, a random variable $X$ is called subgaussian, 
if $\E\,e^{cX^2} < \infty$ for some constant $c > 0$. Assuming that $X$ 
has mean zero, this is equivalent to the statement that the moment 
generating function (or the two-sided Laplace transform) of $X$ satisfies
\be
\E\,e^{tX} \leq e^{\sigma^2 t^2/2}, \quad t \in \R,
\en
with some constant $\sigma^2$. Its optimal value appears in the literature
under different names such as a subgaussian constant or as an optimal proxy variance.
Being deepely connected with logarithmic Sobolev constants and concentration
of measure phenomena, the problem of computation or estimation of $\sigma$
is of considerable interest (including a similar quantity for 
a more general setting of metric spaces, cf. \cite{B-G}).

For example, in the case of a centered Bernoulli distribution $p\delta_q + q \delta_{-p}$,
the subgaussian constant was identified, although not with a rigorous proof,
by Kearns and Saul \cite{K-S} to be
\be
\sigma^2 = \frac{p-q}{2\,(\log p - \log q)}
\en
(cf. also \cite{B-H-T}, \cite{B-K}).
A similar expression was obtained by Diaconis and Saloff-Coste \cite{D-SC} and
by Higuchi and Yoshida \cite{H-Y} for the logarithmic Sobolev constant
of a Markov chain on the two point space. 

Immediate consequences of inequality (1.1) are the finiteness of moments 
of all orders of $X$ and in particular the relations $\E X = 0$ and $\E X^2 \leq \sigma^2$, 
which follow by an expansion of both sides of (1.1) around $t=0$. 
Here the possible case  $\sigma^2 = \Var(X)$ is of particular interest.
The following definition seemed to have appeared first in the work of
Buldygin and Kozachenko \cite{Bu-K1} who called this property ``strongly subgaussian".

\vskip5mm
{\bf Definition.} The random variable $X$ is called {\it strictly
subgaussian}, or the distribution of $X$ is {\it strictly subgaussian}, 
if (1.1) holds with the optimal constant $\sigma^2 = \Var(X)$.

\vskip5mm
Such distributions appear in a natural way in a variety of mathematical problems, 
as well as in statistical mechanics and quantum field theory. For example, under 
the name ``sharp subgaussianity", this class was recently considered
in the work by Guionnet and Husson \cite{G-H} as a condition for LDPs 
for the largest eigenvalue of Wigner matrices with the same rate function as 
in the case of Gaussian entries. Our interest has been motivated by the study
of the central limit theorem with respect to information-theoretic distances.
Let us clarify this connection in the following statement.

Given independent copies $(X_n)_{n \geq 1}$ of a random variable $X$ with 
mean zero and variance one, suppose that the normalized sums
$Z_n = \frac{1}{\sqrt{n}} (X_1 + \dots + X_n)$ have densities $p_n$
for large~$n$. The R\'enyi divergence of order $\alpha > 0$ 
from the distribution of $Z_n$ to the standard normal law with density $\varphi$ 
(or the relative $\alpha$-entropy) is defined by 
\be
D_\alpha(p_n||\varphi) = \frac{1}{\alpha - 1} \log 
\int_{-\infty}^\infty \Big(\frac{p_n(x)}{\varphi(x)}\Big)^\alpha\,
\varphi(x)\, dx.
\en
It is non-decreasing as a function of $\alpha$, representing a strong
distance-like quantity. Here, the case $\alpha=1$ corresponds 
to the relative entropy (Kullback-Leibler's distance) and another important case
$\alpha = 2$ leads to the function of the $\chi^2$-Pearson distance.

\vskip5mm
{\bf Theorem 1.1.} {\sl Suppose that $D_\alpha(p_n||\varphi) <\infty$ for every 
$\alpha$ and some $n = n_\alpha$. For the convergence 
$D_\alpha(p_n||\varphi) \rightarrow 0$ as $n \rightarrow \infty$
with an arbitrary $\alpha>0$, it is necessary and sufficient that
$X$ is strictly subgaussian.
}

\vskip5mm
This characterization follows from the results of \cite{B-C-G1}, which will be
discussed later.

Of course, the property of being strictly subgaussian does not require
that the distribution of $X$ has a density. From (1.2) it already follows
that the symmeric Bernoulli distribution belongs to this class.
More examples are discussed in Arbel, Marchal and Nguyen \cite{A-M-N},
where it is also shown that the distribution of $X$ does not need be symmetric.
The problem of characterization of the whole class of such distributions
is still open and seems to be highly non-trivial. Nevertheless, there is
a simple general sufficient condition for the strict subgaussianity given by 
Newman \cite{N1} (Theorem 4, see also \cite{Bu-K2}, Chapter 1)
in terms of location of zeros of the characteristic function
$$
f(z) = \E\,e^{izX}, \quad z \in \C.
$$
Note that the subgaussian property (1.1) ensures that $f$ has an analytic
extension from the real line to the whole complex plane as an entire function
of order at most 2.

\vskip5mm
{\bf Theorem 1.2.} {\sl Let $X$ be a subgaussian random variable with 
mean zero. If all zeros of $f(z)$ are real, then $X$ is strictly subgaussian. 
}

\vskip5mm
This condition is easily verified for many interesting classes including, 
for example, arbitrary Bernoulli sums and (finite or infinite) convolutions 
of uniform distributions on bounded symmetric intervals. 

The probability distributions of Theorem 1.2 form an important class 
$\mathfrak L$, introduced and studied by Newman in the mid 1970's in 
connection with the Lee-Yang property which naturally arises in the context 
of ferromagnetic Ising models, cf. \cite{N1,N2,N3,N-W}. We will recall
the argument and several properties of this class in Section 4.

Note that if the characteristic function $f(z)$ of a subgaussian random
variable $X$ does not have any real or complex zeros, a well-known theorem 
due to Marcinkiewicz \cite{M} implies that the distribution of $X$ is already Gaussian.
Thus, non-normal subgaussion distributions need to have zeros. 
Towards the characterization problem, the main purpose of this note is 
to explore two natural subclasses of distributions outside $\mathfrak L$ 
that are still strictly subgaussian.
First, we extend Theorem 1.2 in terms of zeros of characteristic functions.

\vskip5mm
{\bf Theorem 1.3.} {\sl Let $X$ be a subgaussian random variable with 
symmetric distribution. If all zeros of $f(z)$ with ${\rm Re}(z) \geq 0$ lie in the
cone centered on the real axis defined by
\be
|{\rm Arg}(z)| \leq \frac{\pi}{8},
\en 
 then $X$ is strictly subgaussian. 
}

\vskip5mm
At the first sight, the condition (1.4) looks artificial. However, it turns out to be
necessary in the following simple situation:

\vskip5mm
{\bf Theorem 1.4.} {\sl Let $X$ be a random variable with a symmetric 
subgaussian distribution. Suppose that $f$ has exactly one zero $z = x+iy$
in the positive quadrant $x,y \geq 0$. Then $X$ is strictly subgaussian, 
if and only if $(1.4)$ holds true.
}

\vskip5mm
As a consequence of Theorem 1.3,
one can partially address the following question from the theory of entire 
characteristic functions (which is one of the central problems in this area): 
What can one say about the possible location of zeros of such functions?

\vskip5mm
{\bf Theorem 1.5.} {\sl Let $(z_n)$ be a finite or 
infinite sequence of non-zero complex numbers in the angle
$|{\rm Arg}(z_n)| \leq \frac{\pi}{8}$ such that
$$
\sum_n \frac{1}{|z_n|^2} < \infty.
$$
Then there exists a symmetric strictly subgaussian distribution
whose characteristic function has zeros exactly at the points
$\pm z_n$, $\pm \bar z_n$. 
}

\vskip5mm
It will be shown that a random variable $X$ with such distribution may be 
constructed as the sum $X = \sum_n X_n$ of independent strictly 
subgaussian random variables $X_n$ whose characteristic function
has zeros at the points $\pm z_n$, $\pm \bar z_n$ (and only 
at these points like in Theorem 1.4). Moreover, one may require that
$$
\Var(X) = \Lambda \sum_n \frac{1}{|z_n|^2}
$$
with any prescribed value $\Lambda \geq \Lambda_0$ where
$\Lambda_0$ is a universal constant ($\Lambda_0 \sim 5.83$).

Returning to Theorem 1.3, it will actually be shown that, 
if a strictly subgaussian random variable $X$ is not normal, the inequality (1.1) 
may be further sharpened as follows: For any $t_0 > 0$, there exists $c = c(t_0)$,
$0 < c < \sigma^2 = \Var(X)$, such that
\be
\E\,e^{tX} \leq e^{c t^2/2}, \quad |t| \geq t_0.
\en
In particular, such a refinement applies to Theorem 1.2. The property (1.5)
is important in the study of rates in the
local limit theorems such as CLT for the R\'enyi divergence of infinite order.
Two results in this direction will be mentioned in the end of this note.

The sharpening (1.5) raises the question of whether or not this separation-type 
property is fulfilled automatically for any non-normal strictly subgaussian 
distribution. At least, it looks natural to expect  the weaker relation
\be
L(t) < e^{\sigma^2 t^2/2}, \quad t \neq 0,
\en
for the Laplace transform $L(t)=\E \exp\{tX\}$. However, the answer to this
is negative, and moreover, (1.6) may turn into an equality for infinitely many 
points $t$. In addition, the characteristic function $f(z)$ 
may have infinitely many zeros approaching the imaginary line
${\rm Arg}(z) = \frac{\pi}{2}$. To this aim, we introduce the following:

\vskip5mm
{\bf Definition.} We say that the distribution $\mu$ of a random variable $X$ 
is periodic with respect to the standard normal law $\gamma$, with period 
$h>0$, if it has a density $p(x)$ such that the density of $\mu$
with respect to $\gamma$,
$$
q(x) = \frac{p(x)}{\varphi(x)} = \frac{d\mu(x)}{d\gamma(x)},
\quad x \in \R,
$$
represents a periodic function with period $h$, that is,
$q(x+h) = q(x)$ for all $x \in \R$.

\vskip5mm
We denote the class of all such distributions 
by $\mathfrak F_h$ and say that $X$ belongs to $\mathfrak F_h$.
The following characterization in terms of Laplace transforms may be useful.

\vskip5mm
{\bf Theorem 1.6.} {\sl Any random variable $X$ in $\mathfrak F_h$
is subgaussian, and the Laplace transform of its distribution is resepresentable as
\be
L(t) = \Psi(t)\,e^{t^2/2}, \quad t \in \R,
\en
where the function $\Psi$ is periodic with period $h$. Conversely, if $\Psi(t)$ 
for a subgaussian random variable $X$ is $h$-periodic, then $X$ belongs 
to $\mathfrak F_h$, as long as the characteristic function $f(t)$ of $X$ is integrable.
}

\vskip5mm
In this way, we obtain a wide class of strictly subgaussian distributions, by requiring
that $\Psi(t) \leq 1$ for all $t$. As a simple example, for any sufficiently small $c>0$,
$$
L(t) = (1 - c \sin^4 t)\,e^{t^2/2}, \quad f(z) = (1 - c \sinh^4 z)\,e^{-z^2/2},
$$
represent respectively the Laplace transform and the characteristic function of 
a strictly subgaussian distribution with mean zero and variance one. In this case, 
we have $L(t) = e^{t^2/2}$ for all $t = \pi k$, $k \in \Z$, and
$f(z_m) = 0$ for $z_m = a + 2\pi i m$, $m \in \Z$, where $a>0$ depends on 
the parameter $c$. Hence ${\rm Arg}(z_m) \rightarrow \frac{\pi}{2}$
as $m \rightarrow \infty$.

More examples based on trigonometric polynomials will be described in Section 12.
The proof of Theorem 1.6 is given in  Sections 10-11. Theorems 1.3 and 1.5 are proved
in Sections 8-9, with preliminary steps in Sections 6-7, and Section 5 is devoted to the
proof of Theorem 1.4. In Sections 3-4 we recall basic definitions and results related 
to the Hadamard factorization theorem the class $\mathfrak L$. 
We conclude with some remarks on the central limit theorem with respect 
to the R\'enyi divergences. Thus, our plan is the following:

\vskip2mm
1. Introduction

2. Basic properties and examples of strictly subgaussian distributions

3. Hadamard's and Goldberg-Ostrovskiĭ's theorems

4. Characteristic functions with real zeros

5. More examples of strictly subgaussian distributions

6. One characterization of characteristic functions

7. Strictly subgaussian symmetric distributions with

\hskip5mm characteristic functions having exactly one non-trivial zero

8. General case of zeros in the angle $|{\rm Arg}(z)| \leq \frac{\pi}{8}$

9. Proof of Theorem 1.5

10. Laplace transforms with periodic components

11. Proof of Theorem 1.6

12. Examples involving triginometric series

13. Examples involving Poisson formula and theta functions

14. Central limit theorems for R\'enyi distances

\vskip5mm
\section{{\bf Basic Properties and Examples of Strictly Subgaussian Distributions}}
\setcounter{equation}{0}

\vskip2mm
\noindent
In addition to the properties $\E X = 0$ and $\E X^2 \leq \sigma^2$, 
the Taylor expansion of the exponential function in (1.1) around zero 
implies as well  that necessarily $\E X^3 = 0$ and $\E X^4 \leq 3\sigma^4$.
Here an equality is attained for symmetric normal distributions
(but not exclusively so).

Turning to other properties and some examples, first let us emphasize
the following two immediate consequences of (1.1).

\vskip5mm
{\bf Proposition 2.1.} {\sl If the random variables $X_1,\dots,X_n$ are 
independent and strictly subgaussian, then their sum $X = X_1 + \dots + X_n$ 
is strictly subgaussian, as well.
}

\vskip5mm
{\bf Proposition 2.2.} {\sl If strictly subgaussian random variables 
$(X_n)_{n \geq 1}$ converge weakly in distribution to a random variable $X$ with
finite second moment, and $\Var(X_n) \rightarrow \Var(X)$ as $n \rightarrow \infty$,
then $X$ is strictly subgaussian.
}

\vskip5mm
{\bf Proof.} By the assumption, putting $\sigma_n^2 = \Var(X_n)$, we have
\be
\E\,e^{tX_n} \leq e^{\sigma_n^2 t^2/2}, \quad t \in \R.
\en
By the weak convergence, $\lim_{n \rightarrow \infty} \E\, u(X_n) = \E\, u(X)$
for any bounded, continuous function $u$ on the real line. 
In particular, for any $c \in \R$,
$$
\lim_{n \rightarrow \infty} \E\,e^{t \min(X_n,c)} = \E\,e^{t \min(X,c)}. 
$$
Hence, by (2.1), for any $t \in \R$,
$$
\E\,e^{t \min(X,c)} \leq \liminf_{n \rightarrow \infty}\, 
\E\,e^{tX_n} \leq \liminf_{n \rightarrow \infty}\, e^{\sigma_n^2 t^2/2} = e^{\sigma^2 t^2/2},
$$
where $\sigma^2 = \Var(X)$. Letting $c \rightarrow \infty$, we get (1.1).
\qed

\vskip2mm
Combining Proposition 2.1 with Proposition 2.2, we obtain:

\vskip5mm
{\bf Corollary 2.3.} {\sl If $\sum_{n=1}^\infty \Var(X_n) < \infty$ for independent, 
strictly subgaussian summands $X_n$, then the series
$X = \sum_{n=1}^\infty X_n$ represents a strictly subgaussian random variable.
}

\vskip5mm
Here, the variance assumption
ensures that the series $\sum_{n=1}^\infty X_n$ is convergent with
probability one (by the Kolmogorov theorem), so that the partial sums
of the series are weakly convergent to the distribution of $X$.
Thus, the class of strictly subgaussian distributions is closed in the weak
topology under infinite convolutions. 

Obviously, it is also closed when taking convex mixtures.

\vskip5mm
{\bf Proposition 2.4.} {\sl If $X_n$ are strictly subgaussian random 
variables with $\Var(X_n) = \sigma^2$, and $\mu_n$ are distributions 
of $X_n$, then for any sequence $p_n \geq 0$ such that 
$\sum_{n=1}^\infty p_n = 1$, the random variable with distribution
$\mu = \sum_{n=1}^\infty p_n \mu_n$
is strictly subgaussian as well and has variance $\Var(X) = \sigma^2$.
}

\vskip5mm
Note also that, if $X$ is strictly subgaussian, then
$\lambda X$ is strictly subgaussian for any $\lambda \in \R$.

Finally, let us give a simple sufficient condition for the property (1.5).
Recall the notation $K(t) = \log \E\,e^{tX}$, $t \in \R$.

\vskip5mm
{\bf Proposition 2.5.} {\sl Let $X$ be a non-normal strictly subgaussian random 
variable. If the function $K(\sqrt{|t|})$ is concave on the half-axis 
$t>0$ and is concave on the half-axis $t<0$, then $(1.5)$ holds true.
}

\vskip5mm
{\bf Proof.} Let $\Var(X) = \sigma^2$. For $t \geq 0$, write
$$
\E\,e^{tX} = e^{\frac{1}{2}\,\sigma^2 t^2 - W(t^2)}.
$$
By the assumption, $W(s)$ is non-negative and convex
in $s\geq 0$, with $W(0) = 0$. In addition, it is $C^\infty$-smooth 
on $(0,\infty)$. Since $X$ is not normal, necessarily $W(s)>0$ and
$W'(s)>0$ for all $s>0$. Using that $W'(s) \uparrow r$ as 
$s \rightarrow \infty$ for some $r \in (0,\infty]$, it follows that
$$
r(s) \, \equiv \, \frac{1}{s} W(s) \, = \,
\int_0^1 W'(sv)\,dv \, \uparrow \, r \quad {\rm as} \ s \rightarrow \infty.
$$
In particular, given $s_0 > 0$, we have
$\frac{1}{s} W(s)  \geq r(s_0) > 0$ for all $s \geq s_0$, or equivalently
$$
K(t) \leq \Big(\frac{1}{2} \sigma^2 - r(s_0)\Big) t^2, \quad t \geq \sqrt{s_0},
$$
which is the desired conclusion. A similar argument works for
$t<0$ as well.
\qed

\vskip2mm
An application of Corollary 2.3 allows to construct a rather rich family
of probability distributions from the class $\mathfrak L$. 
Recall that $L(t) = \E\,e^{tX}$ denotes the Laplace transform.

\vskip5mm
{\bf Example 2.6.} First of all, if a random variable $X$ has a normal distribution 
with mean zero and variance $\sigma^2$,
then it is strictly subgaussian with $L(t) = e^{\sigma^2 t^2/2}$, $t \in \R$.

\vskip2mm
{\bf Example 2.7.} If $X$ has a symmetric Bernoulli distribution, supported on two
points $\pm 1$, then it is strictly subgaussian with
$L(t) = \cosh(t) = \frac{e^t + e^{-t}}{2}$.

\vskip2mm
{\bf Example 2.8.}  If $X = \sum_{n=1}^\infty a_n X_n$ is a Bernoulli sum, 
$\P\{X_n = \pm 1\} = \frac{1}{2}$, $\sum_{n=1}^\infty a_n^2 < \infty$,
with $X_n$ independent, then it is strictly subgaussian with variance 
$\sigma^2 = \sum_{n=1}^\infty a_n^2$. The Laplace transform and 
characteristic function $f$ of $X$ are  given by 
$$
L(t) = \prod_{n=1}^\infty \cosh(a_n t), \quad 
f(t) = \prod_{n=1}^\infty \cos(a_n t).
$$

{\bf Example 2.9.} If $X$ is uniformly distributed on an
interval $[-a,a]$, $a>0$, it is strictly subgaussian. In this case
it may be represented (in the sense of distributions) as the sum 
$$
X = \sum_{n=1}^\infty \frac{a}{2^n}\, X_n, \quad 
\P\{X_n = \pm 1\} = \frac{1}{2} \quad (X_n \ {\rm independent}).
$$
Hence, this case is covered by the previous example, with $L(t) = \frac{\sinh(at)}{at}$.

{\bf Example 2.10.}  If the random variables $X_n$ are independent and 
uniformly distributed on the interval $[-1,1]$, then the infinite sum
$X = \sum_{n=1}^\infty a_n X_n$ with $\sum_{n=1}^\infty a_n^2 < \infty$
represents a strictly subgaussian random variable with
$$
L(t) = \prod_{n=1}^\infty \frac{\sinh(a_n t)}{a_n t}.
$$

{\bf Example 2.11.}  Suppose that $X$ has density $p(x) = x^2 \varphi(x)$, where
$\varphi(x) = \frac{1}{\sqrt{2\pi}}\,e^{-x^2/2}$ is the standard normal
density. Then $\E X = 0$, $\sigma^2 = \E X^2 = 3$, and
$$
L(t) = (1 + t^2)\,e^{t^2/2} \leq e^{3t^2/2}.
$$
Hence, $X$ is strongly subgaussian.

{\bf Example 2.12.}  More generally, if $X$ has a  density of the form
$$
p(x) = \frac{1}{(2d-1)!!}\,x^{2d} \varphi(x), \quad x \in \R, \ d = 1,2,\dots,
$$ 
then 
$\E X = 0$, $\sigma^2 = \E X^2 = 2d+1$, and
$$
L(t) = \frac{1}{(2d-1)!!}\,H_{2d}(it)\,e^{t^2/2} \leq 
e^{(2d+1)\,t^2/2}.
$$
Hence, $X$ is strictly subgaussian. The last inequality follows from
Theorem 1.2, since the Chebyshev-Hermite polynomials have real zeros, only.

\vskip5mm
\section{{\bf Hadamard's and Goldberg-Ostrovskiĭ's Theorems}}
\setcounter{equation}{0}

\vskip2mm
\noindent
All the previous examples may be included as partial cases of a more general setup.
First, let us recall some basic definitions and notations related to the Hadamard 
theorem from the theory of complex variables. Given an entire function $f(z)$, introduce
$$
M_f(r) = \max_{|z| \leq r} |f(z)| = \max_{|z| = r} |f(z)|, \quad r \geq 0,
$$
which characterizes the growth of $f$ at infinity. The order of $f$ is defined by
$$
\rho = \limsup_{r \rightarrow \infty} \frac{\log \log M_f(r)}{\log r}.
$$
Thus, $\rho$ is an optimal value such that, for any $\ep > 0$, we have
$M_f(r) < e^{r^{\rho + \ep}}$ for all large $r$. If $f$ is a polynomial, then 
$\rho = 0$. If $\rho$ is finite, then the type of $f$ is defined by
$$
\tau = \limsup_{r \rightarrow \infty} \frac{\log M_f(r)}{r^\rho}.
$$
Thus, $\tau$ is an optimal value such that, for any $\ep > 0$, we have
$M_f(r) < e^{(\tau + \ep)\,r^\rho}$ for all sufficiently large $r$.
If $0 < \tau < \infty$, the function $f$ is said to be of normal type.

For integers $p \geq 0$, introduce the functions
$$
G_p(u) = (1 - u)\exp\Big\{u + \frac{u^2}{2} + \dots + \frac{u^p}{p}\Big\},
\quad u \in \C,
$$
called the primary factors, with the convention that $G_0(u) = 1-u$.
Given a sequence of complex numbers $z_n \neq 0$ such that 
$|z_n| \uparrow \infty$, one considers a function of the form 
\be
\Pi(z) = \prod_{n = 1}^\infty G_p(z/z_n)
\en
called a canonical product. An integer $p \geq 0$ is called 
the genus of this product, if it is the smallest integer such that
\be
\sum_{n=1}^\infty \frac{1}{|z_n|^{p+1}} < \infty.
\en

There is a simple estimate $\log |G_p(u)| \leq A_p |u|^{p+1}$ 
where the constant $A_p$ depends on $p$ only. Therefore, the product
in (3.1) is uniformly convergent as long as (3.2) is fulfilled.

See e.g. Levin \cite{Le} for the following classical theorem.

\vskip5mm
{\bf Theorem 3.1} (Hadamard). {\sl Any entire function $f$ of a finite order $\rho$
can be represented in the form
\be
f(z) = z^m\,e^{P(z)} \prod_{n \geq 1} G_p(z/z_n), \quad z \in \C.
\en
Here $z_n$ are the non zero roots of $f(z)$, the genus of the canonical
product satisfies $p \leq \rho$, $P(z)$ is a polynomial of degree 
$\leq \rho$, and $m \geq 0$ is the multiplicity of the zero at the origin.
}

\vskip5mm
In order to describe the convergence of the canonical product, 
assume that $f(z)$ has an infinite sequence of non-zero
roots $z_n$ arranged in increasing order of their moduli so that
$$
0 < |z_1| \leq |z_2| \leq \dots \leq |z_n| \rightarrow \infty \quad 
{\rm as} \ n \rightarrow \infty.
$$
Define the convergence exponent of the sequence $a_n$ by
$$
\rho_1 = \inf\Big\{\lambda > 0:
\sum_{n=1}^\infty \frac{1}{|z_n|^\lambda} < \infty\Big\}.
$$
A theorem due to Borel asserts that the order $\rho$ of the canonical product $\Pi(z)$
satisfies $\rho \leq \rho_1$. Moreover, Theorem 6 from \cite{Le}, p.16, 
states that the convergence exponent of the zeros of any entire function 
$f(z)$ does not exceed its order: $\rho_1 \leq \rho$. Thus,
for canonical products the convergence exponent of the zeros is equal to the
order of the function: $\rho_1 = \rho$ (Theorem 7).

There is also the following elementary relation between the convergence exponent
and the genus of the canonical product: $p \leq \rho_1 \leq p+1$.
Assuming that $\rho_1$ is an integer, we have that
$\sum_{n=1}^\infty |z_n|^{-\rho_1} = \infty \Rightarrow p = \rho_1$, 
while $p = \rho_1 + 1$ means that the latter series is convergent.

The following theorem due to Goldberg and Ostrovskiĭ \cite{G-O} 
refines Theorem 3.1 for the class of ridge entire functions 
whose all zeros are real. Recall that $f$ is a ridge function, if it satisfies
$|f(x+iy)| \leq |f(iy)|$ for all $x,y \in \R$.

\vskip5mm
{\bf Theorem 3.2} (Goldberg-Ostrovskiĭ). {\sl Suppose that an entire ridge function
$f$ of a finite order has only real roots. Then it can be represented in the form
\be
f(z) = c\,e^{i\beta z - \gamma z^2/2} \prod_{n \geq 1} 
\Big(1 - \frac{z^2}{z_n^2}\Big), \quad z \in \C,
\en
for some $c \in \C$, $\beta \in \R$, $\gamma \geq 0$, and $z_n > 0$
such that $\sum_{n \geq 1} z_n^{-2} < \infty$.
}

\vskip5mm
We refer to \cite{G-O}. See also Kamynin \cite{K} for generalizations of 
Theorem 3.2 to the case where the zeros of $f$ are not necessarily real.

\vskip5mm
\section{{\bf Characteristic Functions with Real Zeros}}
\setcounter{equation}{0}

\vskip2mm
\noindent
We are now prepared to prove Theorem 1.2, including the relation
(1.5) in the non-Gaussian case which is stronger than (1.1).

Thus, let $X$ be a subgaussian random variable with mean zero
and variance $\sigma^2 = \Var(X)$. Then the inequality (1.1)
may be extended to the complex plane in the form
$$
|f(z)| \leq e^{b\, {\rm Im}(z)^2/2}, \quad z \in \C,
$$
for some constant $b \geq \sigma^2$, where $f$ is the characteristic function
of $X$. Hence, $f$ is a ridge entire function of order $\rho \leq 2$. We are 
therefore able to apply Theorem 3.2 which yields the representation (3.4)
for some $c \in \C$, $\gamma \geq 0$, $\beta \in \R$, and for some
finite or infinite sequence $z_n > 0$ such that 
$\sum_{n \geq 1} z_n^{-2} < \infty$. Note that
$f(z_n) = f(-z_n) = 0$, so that $\{z_n,-z_n\}$ are all zero of $f$
(this set may be empty).
Since $f(0) = 1$ and $f'(0) = 0$, we necessarily have $c = 1$ and
$\beta = 0$. Hence, this representation is simplified to
\be
f(z) = e^{-\gamma z^2/2} \prod_{n \geq 1} 
\Big(1 - \frac{z^2}{z_n^2}\Big).
\en
Since $f''(0) = -\sigma^2$, we also have
\be
\frac{1}{2} \sigma^2 = \frac{1}{2} \gamma + \sum_{n \geq 1} \frac{1}{z_n^2},
\en
so that $\gamma \leq \sigma^2$. Applying (4.1) with 
$z = -it$, $t \in \R$, we get a similar representation for 
the Laplace transform
\be
\E\,e^{tX} = e^{\gamma t^2/2} \prod_{n \geq 1} 
\Big(1 + \frac{t^2}{z_n^2}\Big).
\en
Using $1+x \leq e^x$ ($x \in \R$), we see that the right-hand side 
above does not exceed $e^{\sigma^2 t^2/2}$,
where we used (4.2). Hence (4.3) leads to the desired bound (1.1), 
and Theorem 1.2 is  proved.

Let us also verify the property (1.5) in the case where the random
variable $X$ is not normal. Then the product in (4.3) is not empty 
and therefore $\gamma < \sigma^2$. Let us rewrite (4.3) as
$$
\E\,e^{tX} = e^{V(t^2)}, \quad V(s) = \gamma s + \sum_{n \geq 1} 
\log \Big(1 + \frac{s}{z_n^2}\Big).
$$
Since the function $V$ is concave, it remains to refer to Proposition 2.5.
\qed

\vskip2mm
{\bf Remark.} Using (4.2), let us rewrite (4.1) with $z = t \in \R$ 
in the form
\be
f(t) = e^{-(3 \gamma - \sigma^2)\, t^2/4} \prod_{n \geq 1} 
\Big(1 - \frac{t^2}{z_n^2}\Big)\,e^{-\frac{t^2}{2z_n^2}}.
\en
Here, the terms in the product represent characteristic funtions of
random variables $\frac{1}{z_n} X_n$ such that all $X_n$ have
density $p(x) = x^2 \varphi(x)$ which we discussed in Example 2.11.
Hence, if $\frac{1}{3} \sigma^2 \leq \gamma \leq \sigma^2$,
the function $f(t)$ in (4.4) represents the characteristic function of 
$$
X = cZ + \sum_{n \geq 1} \frac{1}{z_n} X_n, \quad 
c^2 = \frac{3}{2} \gamma - \frac{1}{2}\sigma^2,
$$
assuming that $X_n$ are independent and $Z \sim N(0,1)$ is 
independent of all $X_n$.

Note that (4.1) does not always define a characteristic function. 
For example, when there is only one term in the product, we have
$f(t) = e^{-\gamma t^2/2} (1 - \frac{t^2}{z_1^2})$.
It is a characteristic function, if and only if 
$\gamma \geq \frac{1}{z_1^2}$ (cf. e.g. \cite{L-O}, p.\,34).
We will return to this question in Section 8.

\vskip2mm
{\bf Properties of the class} $\mathfrak L$. Following Newman \cite{N1}, 
let us emphasize several remarkable properties of strictly subgaussian
distributions whose characteristic functions have real zeros, only. 
Starting from (4.3), one can represent the log-Laplace transform of $X$ as
$$ 
K(t) = \log \E\,e^{tX} = \frac{\gamma t^2}{2} + \sum_{n \geq 1} 
\log \Big(1 + \frac{t^2}{z_n^2}\Big).
$$
Hence the cumulants of even order $2m$ of $X$ are given for $m \geq 2$ 
by
$$
\gamma_{2m} = K^{(2m)}(0) = (-1)^{m-1}\,\frac{(2m)!}{m}
\sum_{n \geq 1} \frac{1}{z_n^{2m}},
$$
while $\gamma_{2m-1} = 0$. In particular, the distribution of $X$
has to be symmetric about the origin, with $(-1)^{m-1} \gamma_{2m} \geq 0$.
As was also shown in \cite{N1}, the cumulants and the moments 
of $X$ admit the bounds
\be
(-1)^{m-1} \gamma_{2m} \leq \frac{(2m)!}{2^m \, m}\,\sigma^{2m}, \quad
\E X^{2m} \leq \frac{(2m)!}{2^m \, m!}\,\sigma^{2m}.
\en
In addition, for all integers $k \geq 0$ and $t \in \R$,
$$
\sum_{m=1}^{2k} \frac{\gamma_{2m}}{(2m)!}\, t^{2m} \leq K(t) \leq 
\sum_{m=1}^{2k+2} \frac{\gamma_{2m}}{(2m)!}\, t^{2m}.
$$

Since the class $\mathfrak L$ is closed under convolutions, the second
inequality in (4.5) continues to hold for weighted sums of independent, 
strictly subgaussian random variables. This provides a natural extension
of Khinchine's inequality for Bernoulli sums, as noticed in \cite{N2}
(cf. also a recent work \cite{H-N-T}).

\vskip5mm
\section{{\bf More Examples of Strictly Subgaussian Distributions}}
\setcounter{equation}{0}

\vskip2mm
\noindent
In connection with the problem of location of zeros,
we now examine probability distributions with characteristic
functions of the form
\be
f(t) = e^{-t^2/2}\,(1 - \alpha t^2 + \beta t^4),
\en
where $\alpha,\beta \in \R$ are parameters. It was already mentioned
that when $\beta = 0$, we obtain a characteristic function 
$$
f(t) = e^{-t^2/2}\,(1 - \alpha t^2),
$$
if and only if $0 \leq \alpha \leq 1$.  As we will see, in the general case,
it is necessary that $\beta \geq 0$ for $f(t)$ to be a characteristic function
(although negative values of $\alpha$ are possible for small $\beta$).
Before deriving a full characterization, first let us emphasize the following.

\vskip5mm
{\bf Proposition 5.1.} {\sl Given $\beta \geq 0$, a random variable $X$
with characteristic function of the form $(5.1)$ is strongly subgaussian, 
if and only if $\alpha$ satisfies $\alpha \geq \sqrt{2\beta}$.
}

\vskip5mm
{\bf Proof.} Recall that $X$ is strongly subgaussian, if and only if,
for all $t \in \R$,
\be
\E\,e^{tX} \leq e^{\sigma^2 t^2/2}, \quad \sigma^2 = -f''(0).
\en
Near zero, the characteristic function in (5.1) behaves 
like a quadratic polynomial
$
f(t) = 1 - \frac{1}{2}\,t^2 - \alpha t^2 + O(t^4),
$
so that $\sigma^2 = 1 + 2\alpha$ (in particular, $\alpha \geq -\frac{1}{2}$). 
Hence, applying (5.1) to the values $-it$, one may rewrite (5.2) equivalently
(multiplying both sides by $\exp(-t^2/2)$) as
$$
1 + \alpha t^2 + \beta t^4 \leq e^{\alpha t^2} = 
1 + \alpha t^2 + \frac{1}{2}\,\alpha^2 t^4 + 
\frac{1}{6}\,\alpha^3 t^6 + \dots
$$
If $\alpha \geq 0$, this inequality holds for all $t \in \R$, if and 
only if $\alpha^2 \geq 2\beta$. As for the case $\alpha < 0$, this is
impossible, since then $e^{\alpha t^2} \rightarrow 0$ as
$t \rightarrow \infty$ exponentially fast.
\qed

\vskip5mm
As already emphasized, if a random variable $X$ is subgaussian 
(even if it is not strictly subgaussian), its characteristic function 
$f(t)$ may be extended to the complex plane as an entire function
$f(z) = \E\,e^{izX}$ of order $\rho \leq 2$ and of finite type 
like in the strictly subgaussian case (5.2). Since in general
$f(-\bar z) = \bar f(z)$, any zero $z = x+iy$ of $f$ ($x,y \in \R$) is 
complemented with zero $-\bar z = -y-ix$. If in addition the distribution
of $X$ is symmetric about zero, then $-z$ and $\bar z$ will also be
zeros of $f$. Thus, in this case with every non-real zero $z$, 
the characteristic function has 3 more distinct zeros, and hence we have
4 distinct zeros $\pm x \pm iy$, $x,y > 0$. One can now apply
Proposition 5.1 to prove Theorem 1.4.

\vskip5mm
{\bf Proof of Theorem 1.4.} Given a random variable $X$ with 
a symmetric subgaussian distribution, suppose that its characteristic
function has exactly one zero $z = x+iy$ in the positive quadrant
$x,y \geq 0$. We need to show that $X$ is strongly subgaussian, 
if and only if
\be
0 \leq {\rm Arg}(z) \leq \frac{\pi}{8}.
\en

The case where $z = x$ is real is covered by Theorem 1.2. 
The argument below also works in this case, but for definiteness let us assume 
that $z$ is complex, so that $x,y>0$ (the case $x = 0$ and $y>0$ is 
impossible, since then $f(z) = f(iy) \geq 1$).

Thus, let $f(z)$ have four distinct roots $z_1 = z$, 
$z_2 = -z = -x - iy$, $z_3 = \bar z = x-iy$, $z_4 = - \bar z = -x+iy$.
Applying Hadamard's theorem, we get a representation
$$
f(z) = e^{P(z)}\,\Big(1 - \frac{z}{z_1}\Big)\Big(1 - \frac{z}{z_2}\Big)
\Big(1 - \frac{z}{z_3}\Big)\Big(1 - \frac{z}{z_4}\Big),
$$
where $P(z)$ is a quadratic polynomial. Since $f(0)=1$, 
necessarily $P(0) = 0$. Also, by the symmetry of the distribution 
of $X$, we have $f(z) = f(-z)$, which implies $P(z) = P(-z)$ for all
$z \in \C$. It follows that $P(z)$ has no linear term, so that
$P(z) = -\frac{1}{2}\gamma z^2$ for some $\gamma \in \C$. Thus, putting 
$w = a+bi = \frac{1}{x+iy}$, we have
\begin{eqnarray}
f(t) 
 & = &
e^{-\gamma t^2/2}\,(1 - w t)(1+w t)(1-\bar w t)(t+\bar w t) \nonumber \\
 & = &
e^{-\gamma t^2/2}\,\big(1 - (w^2 + \bar w^2)\, t^2 + |w|^4 t^4\big) \nonumber \\
 & = &
e^{-\gamma t^2/2}\,\big(1 - 2(a^2 - b^2)\, t^2 + (a^2 + b^2)^2\, t^4\big).
\end{eqnarray}

Comparing both sides of (5.4) near zero according to Taylor's expansion,
we get that 
\be
\gamma + 4(a^2 - b^2) = \sigma^2.
\en
In particular, $\gamma$ must be a real number, necessarily positive
(since otherwise $f(t)$ would not be bounded on the real axis).
Moreover, the case $a=|b|$ is impossible, since then
$f(t) = e^{-\sigma^2 t^2/2}\,(1 + 2b^4 t^4)$.
Rescaling the variable and applying Proposition 5.1 with $\alpha = 0$,
we would conclude that the random variable $X$ is not  strictly subgaussian.

Thus, let $a \neq |b|$ (as we will see, necessarily $\gamma > \sigma^2$). 
Again rescaling of the $t$-variable, one may assume that
$\gamma = 1$ in which case the representation (5.4) becomes
$$
f(t) = e^{-t^2/2}\,\big(1 - 2(a^2 - b^2)\, t^2 + (a^2 + b^2)^2\, t^4\big).
$$
One can now apply Proposition 5.1 with parameters
$\alpha = 2(A-B)$, $\beta = (A+B)^2$, where $A = a^2, \ B = b^2$.
Since the condition $\alpha \geq 0$ is necessary for $f(t)$ to be 
a characteristic function of a strictly subgaussian distribution, 
we may assume that $A \geq B$ (in fact, we have $A>B$, since $a \neq |b|$).
The condition $\beta \leq \frac{1}{2}\,\alpha^2$, that is,
$2 (A-B)^2 \geq (A+B)^2$ is equivalent to 
$$
(A + B)^2 \geq 8 AB \ \Longleftrightarrow \
(a^2 + b^2)^2 \geq 8 a^2 b^2. 
$$
To express this in polar coordinates, put $a = r\cos \theta$, $b = r\sin \theta$
with $r^2 = a^2+b^2$ and $|\theta| \leq \frac{\pi}{2}$. Since 
$A \geq B$, that is $a \geq |b|$, necessarily $|\theta| \leq \frac{\pi}{4}$, 
and the above turns out to be the same as
$$
\cos^2(\theta)\, \sin^2(\theta) \leq \frac{1}{8} \ \Longleftrightarrow \
\sin^2(2\theta) \leq \frac{1}{2} \ \Longleftrightarrow \ |\theta| \leq 
\frac{\pi}{8}.
$$
Since $\theta = {\rm Arg}(a + bi) = -{\rm Arg}(z)$, the desired 
characterization (5.3) follows.
\qed

\vskip5mm
\section{{\bf One Characterization of Characteristic Functions}}
\setcounter{equation}{0}

\vskip2mm
\noindent
It remains to decide whether or not the characteristic functions in 
Proposition 5.1 with non-real zeros do exist. Therefore, we now turn 
to the characterization of the property that the functions of the form
\be
f(t) = e^{-t^2/2}\,(1 - \alpha t^2 + \beta t^4)
\en
are positive definite. The more general class of functions
$f(t) = e^{-\gamma t^2/2}\,(1 - \alpha t^2 + \beta t^4)$, $\gamma>0$,
is reduced to (6.1) by rescaling the $t$-variable.

\vskip5mm
{\bf Proposition 6.1.} {\sl The equality $(6.1)$ defines 
a characteristic function, if and only if the point $(\alpha,\beta)$ 
belongs to one of the following two regions:
\be
4\beta - 2\sqrt{\beta (1 - 2 \beta)} \leq \alpha \leq 3\beta + 1, \qquad 
0 \leq \beta \leq \frac{1}{3}, 
\en
or
\be
4\beta - 2\sqrt{\beta (1 - 2 \beta)} \leq \alpha \leq 
4\beta + 2\sqrt{\beta (1 - 2 \beta)}, \quad 
\frac{1}{3} \leq \beta \leq \frac{1}{2}.
\en
}

\vskip2mm
The expression on the left-hand sides in (6.2)-(6.3) is negative, 
if and only if $\beta < \frac{1}{6}$. Hence, for such values
of $\beta$, the parameter $\alpha$ may be negative.

Combining Propositions 5.1 and 6.1, we obtain a full characterization
of strictly subgaussian distributions with characteristic functions
of the form (6.1). To this aim, one should complement (6.2)-(6.3) 
with the bound $\alpha \geq \sqrt{2\beta}$. To describe the full region, 
we need to solve the corresponding inequalities. First, it should
be clear that $\sqrt{2\beta}$ is smaller than the right-hand sides
of (6.2)-(6.3) for all $0 \leq \beta \leq \frac{1}{2}$. In this 
$\beta$-interval, we also have
$$
4\beta - 2\sqrt{\beta (1 - 2 \beta)} \leq \sqrt{2\beta}
 \, \Longleftrightarrow \, 12\beta - 3 \leq 2\sqrt{2 (1 - 2 \beta)}.
$$
The latter is fulfilled automatically for $\beta \leq \frac{1}{4}$.
For $\frac{1}{4} \leq \beta \leq \frac{1}{2}$, squaring the
above inequality, we arrive at the quadratic  inequality
$$
144 \beta^2 - 56 \beta + 1 \leq 0.
$$
The corresponding quadratic equation has two real roots,
one of which $0.0188...$ is out of our interval, while the other
one
$$
\beta_0 = \frac{1}{36}\,(7 + 2\sqrt{10}) \sim  0.3701...
$$ 
belongs to the interval $(\frac{1}{3},\frac{1}{2})$.
Therefore, the left-hand side in (6.2) should be replaced
with $\sqrt{2\beta}$ on the whole interval 
$0 \leq \beta \leq \frac{1}{3}$, while the lower bounds in (6.3)
should be properly changed for $\beta \leq \beta_0$ and
$\beta \geq \beta_0$. That is, we obtain:

\vskip5mm
{\bf Proposition 6.2.} {\sl The equality $(6.1)$ defines a characteristic 
function of a strictly subgaussian distribution, if and only if
\begin{eqnarray}
\sqrt{2\beta} \leq
 & \hskip-2mm \alpha &  \hskip-2mm  \leq
3\beta + 1, \qquad \qquad \qquad 0 \leq \beta \leq \frac{1}{3}, \\
\sqrt{2\beta}  \leq 
 & \hskip-2mm \alpha &  \hskip-2mm \leq 
4\beta + 2\sqrt{\beta (1 - 2 \beta)}, \quad 
\frac{1}{3} \leq \beta \leq \beta_0, \\ 
4\beta - 2\sqrt{\beta (1 - 2 \beta)} \leq
 & \hskip-2mm \alpha &  \hskip-2mm  \leq
4\beta + 2\sqrt{\beta (1 - 2 \beta)}, \quad 
\beta_0 \leq \beta \leq \frac{1}{2}.
\end{eqnarray}
}

\vskip2mm
{\bf Proof of Proposition 6.1.} Recall that the Chebyshev-Hermite 
polynomial $H_k(x)$ of degree $k = 0,1,2,\dots$ is defined via the identity
$\varphi^{(k)}(x) = (-1)^k H_k(x) \varphi(x)$. In particular, 
$$
H_0(x) = 1, \quad H_2(x) = x^2-1, \quad H_4(x) = x^4 - 6x^2 + 3.
$$
Equivalently, for even orders
$$
t^{2k}\,e^{-t^2/2} = \int_{-\infty}^\infty 
(-1)^k H_{2k}(x) \varphi(x)\,e^{itx}\,dx.
$$
Therefore, the function in (6.1) represents the Fourier transform of the function
\bee
p(x)
 & = &
\big(1 + \alpha H_2(x) + \beta H_4(x)\big)\,\varphi(x) \\
 & = &
\big((1 - \alpha + 3\beta) + (\alpha - 6\beta) x^2 + \beta x^4\big)\,
\varphi(x),
\ene
whose total integral is $f(0) = 1$. Hence, $p(x)$ represents
a probability density, if and only if
$$
\psi(y) \equiv 
(1 - \alpha + 3\beta) + (\alpha - 6\beta) y + \beta y^2 \geq 0 \quad
{\rm for \ all} \ y \geq 0.
$$
Choosing $y = 0$ and $y \rightarrow \infty$, we obtain necessary
conditions
\be
\alpha \leq 3\beta + 1, \quad \beta \geq 0.
\en

Assuming this, a sufficient condition for the inequality $\psi(y) \geq 0$ 
to hold for all $y \geq 0$ is $\alpha \geq 6\beta$. As a result, 
we obtain a natural region for the parameters, namely
\be
6\beta \leq \alpha \leq 3\beta + 1, \quad 0 \leq \beta \leq \frac{1}{3},
\en
for which $f(t)$ in (6.1) is a characteristic function.

In the case $\alpha < 6\beta$, we obtain a second region. 
Note that the quadratic function 
$\psi(y) = c_0 + 2c_1 y + c_2 y^2$ with $c_0,c_2 \geq 0$ and 
$c_1 < 0$ is non-negative in $y \geq 0$, if and only if 
$c_1^2 \leq c_0 c_2$. For the coefficients $c_2 = \beta > 0$
and $2c_1 = \alpha - 6\beta < 0$, the condition $c_1^2 \leq c_0 c_2$
means that
$$
\Big(\frac{\alpha - 6\beta}{2}\Big)^2 \leq (1 - \alpha + 3\beta) \beta
 \, \Longleftrightarrow \, (\alpha - 4 \beta)^2 \leq 4\beta (1 - 2 \beta). 
$$
Thus, necessarily $\beta \leq \frac{1}{2}$, and then admissible
values of $\alpha$ are described by the relations
\be
4\beta - 2\sqrt{\beta (1 - 2 \beta)} \leq \alpha \leq 
4\beta + 2\sqrt{\beta (1 - 2 \beta)}
\en
in addition to the assumption $\alpha < 6\beta$ and the
necessary conditions in (6.7).

If $\frac{1}{3} \leq \beta \leq \frac{1}{2}$, we arrive 
at the desired relations in (6.3), since
$$
4\beta + 2\sqrt{\beta (1 - 2 \beta)} \leq 3\beta + 1 \leq 6\beta.
$$ 

If $\beta \leq \frac{1}{3}$, then $6\beta \leq 3\beta+1$.
In the case $\alpha < 6\beta$, the upper bound in (6.9) 
will hold automatically, since
$$
6\beta \leq 4\beta + 2\sqrt{\beta (1 - 2 \beta)} \quad
{\rm for \ all} \ 0 \leq \beta \leq \frac{1}{3}. 
$$
So, for the values $\alpha < 6\beta$ and $\beta \leq \frac{1}{3}$, 
(6.9) is simplified to
\be
4\beta - 2\sqrt{\beta (1 - 2 \beta)} \leq \alpha \leq 
6\beta, \quad 0 < \beta \leq \frac{1}{3}.
\en
It remains to take the union of the two regions described 
by (6.10) with (6.8), and then we arrive at (6.2).
\qed

\vskip5mm
\section{{\bf Strictly Subgaussian Symmetric Distributions with Characteristic \\
Functions Having Exactly One Non-trivial Zero}}
\setcounter{equation}{0}

\vskip2mm
\noindent
One may illustrate Proposition 6.2 by the following simple example. For 
$\beta = \frac{1}{3}$, admissible values of $\alpha$ cover the interval 
$\sqrt{2/3} \leq \alpha \leq 2$, following both (6.4) and (6.5).
Choosing $\alpha = \sqrt{2/3}$, we obtain the characteristic function 
$$
f(t) = e^{-t^2/2}\,\Big(1 - \sqrt{\frac{2}{3}}\, t^2 + \frac{1}{3} t^4\Big)
$$
of a strictly subgaussian random variable. It has four distinct complex zeros 
$z_k$ defined by $z^2 = r^2\,(1 \pm i)$ with $r^2 = \frac{1}{3}\sqrt{2/3}$, so 
$$
z_1 = (2r)^{1/4}\,e^{i\pi/8}, \quad z_2 = (2r)^{1/4}\,e^{-i\pi/8}, \quad
z_3 = (2r)^{1/4}\,e^{7i\pi/8}, \quad z_4 = (2r)^{1/4}\,e^{-7i\pi/8}.
$$
Note that $|{\rm Arg}(z_{1,2})| = \frac{\pi}{8}$. 
As  already mentioned, it was necessary that 
$|{\rm Arg}(z)| \leq \frac{\pi}{8}$ for all zeros with ${\rm Re}(z) > 0$ 
in the class of all strictly subgaussian 
probability distributions with characteristic functions of the form (6.1).

In order to describe the  possible location of zeros, let us see what 
Proposition 6.2 is telling us about the class of functions
\be
f(t)  = e^{-t^2/2}\,(1 - w t)(1+w t)(1-\bar w t)(t+\bar w t), \quad
t \in \R,
\en
with $w = a+bi$. Thus, in the complex plane 
$f(z)$ has two or four distinct zeros 
$z = \pm 1/w$, $z = \pm 1/\bar w$ depending on whether
$b = 0$ or $b \neq 0$.
Note that 
$$
|{\rm Arg}(z)| = |{\rm Arg}(w)|
$$ 
when $z$ and $w$ are taken from the half-plane 
${\rm Re}(z) > 0$ and ${\rm Re}(w) > 0$.

\vskip5mm
{\bf Proposition 7.1.} {\sl Let $w = a+bi$ with $a > 0$. 
The function $f(t)$ in $(7.1)$ represents a characteristic function 
of a strictly subgaussian random variable, if and only if 
$$
a \leq 2^{-1/4} \sim 0.8409,
$$ 
while 
$|b|$ is sufficiently small. More precisely, this is the case whenever
$|b| \leq b(a)$ with a certain function $b(a) \geq 0$ such that
$b(2^{-1/4}) = 0$ and $b(a)>0$ for $0 < a < 2^{-1/4}$.
}

\vskip2mm
Moreover, there exists a universal 
constant $0 < a_0 < 2^{-1/4}$, $a_0 \sim 0.7391$, such that 
for $0 \leq a \leq a_0$ and only for these $a$-values, the property 
$|b| \leq b(a)$ is equivalent to the angle requirement
${\rm Arg}(w) \leq \frac{\pi}{8}$. As for the values
$a_0 < a \leq 2^{-1/4}$, this angle must be smaller.

\vskip5mm
{\bf Proof.} We may assume that $b \geq 0$.
The function in (7.1) may be expressed in the form
\be
f(t) = e^{-t^2/2}\,(1 - \alpha t^2 + \beta t^4)
\en
with parameters $\alpha = 2(A-B)$, $\beta = (A+B)^2$, where
$A = a^2$, $B = b^2$. Since the condition $\alpha \geq 0$ is 
necessary for $f(t)$ to be a characteristic function of a strictly subgaussian 
distribution, we may require that $a \geq b$, that is, $A \geq B$. 
Recall that
\be
{\rm Arg}(w) \leq \frac{\pi}{8} \, \Longleftrightarrow \,
\alpha \geq \sqrt{2\beta} \, \Longleftrightarrow \,
b \leq \frac{1}{\sqrt{2} + 1}\,a.
\en
In fact, as easy to check, if $w = re^{i\theta}$, then
$$
\alpha^2 - 2\beta = 2\beta \cos(4\theta).
$$

In order to apply Proposition 6.2, first note that the above parameters
satisfy $\alpha \leq 2\sqrt{\beta}$. In this case, the upper bounds
in (6.4)-(6.6) are fulfilled automatically.
Therefore, we only need to take into account the lower bounds
in (6.4)-(6.6). Thus, $f(t)$ in (7.2) represents the characteristic 
function of a strongly subgaussian distribution, if and only if
\be
\frac{1}{\sqrt{2}}\,(A+B) \leq A-B \quad {\rm for} \ \ 
0 < A+B < \sqrt{\beta_0}
\en
or
\be
2(A+B)^2 - (A+B)\sqrt{1 - 2 (A+B)^2} \leq A-B \quad {\rm for} \ \
\sqrt{\beta_0} \leq A+B \leq \frac{1}{\sqrt{2}},
\en 
where $\beta_0 = \frac{1}{36}\,(7 + 2\sqrt{10}) \sim  0.3701...$
Since the condition $A+B \leq 2^{-1/2}$ is necessary,
we should require that $a \leq 2^{-1/4}$. Moreover, for $a = 2^{-1/4}$,
there is only one admissible value $b = 0$,
when $w$ is a real number, $w = 2^{-1/4}$.

Let us recall that
$$
\sqrt{2\beta} < 4\beta - 2\sqrt{\beta(1 - 2\beta)} \quad {\rm for} \ \
\beta_0 < \beta \leq \frac{1}{2},
$$
in which case there is a strict inequality $\alpha > \sqrt{2\beta}$
for admissible values of $\alpha$ in (6.6). Hence, 
${\rm Arg}(w) < \frac{\pi}{8}$ according to (6.3).  Thus, 
${\rm Arg}(w) < \frac{\pi}{8}$ for the region described in (7.5).

Turning to the region of couples $(A,B)$ as in (7.4),
let us fix a value $0 < A < \sqrt{\beta_0}$.
The first inequality in (7.4) is equivalent to
$$
B \leq \frac{\sqrt{2} - 1}{\sqrt{2} + 1}\,A = 
\frac{1}{(\sqrt{2} + 1)^2}\,A,
$$
which is the same as (7.3). The value $B = \frac{1}{(\sqrt{2} + 1)^2}\,A$ 
satisfies the second constraint, if and only if 
$(1 + \frac{1}{(\sqrt{2} + 1)^2})\,A \leq \sqrt{\beta_0}$,
which is equivalent to $0 < a \leq a_0$ with
\be
a_0 = \beta_0^{1/4} \frac{\sqrt{2} + 1}{\sqrt{4 + 2\sqrt{2}}}
\sim 0.7391.
\en
Therefore, in this $a$-interval
Proposition 7.1 holds true with $b(a) = \frac{1}{\sqrt{2} + 1}\,a$.

Now, let $a_0 < a < 2^{-1/4}$. Since $A < \frac{1}{2}$, both (7.4) and (7.5)
are fulfilled for all $B$ small enough. Indeed, if $A \geq \sqrt{\beta_0}$
and $B = 0$, (7.5) becomes
$$
2A^2 - A \sqrt{1 - 2A^2} \leq A,
$$
which holds with a strict inequality sign. To show that (7.5) is
solved as $B \leq B(A)$ for a certain positive function $B(A)$,
it is sufficient to verify that the left-hand side of (7.5) is 
increasing in $B$ (since the right-hand side is decreasing in $B$).
Consider the function
$$
u(x) = 2x^2 - x\sqrt{1 - 2x^2}, \quad 
\sqrt{\beta_0} \leq x < \frac{1}{\sqrt{2}}.
$$
We have
$$
u'(x) = 
4x - \sqrt{1 - 2x^2} + \frac{2x^2}{\sqrt{1 - 2x^2}} \geq 0
$$
for $x \geq \frac{1}{2}$, hence for $x \geq \sqrt{\beta_0}$.
Thus, $u(A+B)$ is increasing in $B$, proving the  claim.
\qed

\vskip5mm
\section{{\bf General Case of Zeros in the Angle 
$|{\rm Arg}(z)| \leq \frac{\pi}{8}$}}
\setcounter{equation}{0}

\vskip2mm
\noindent
We are now prepared to prove Theorem 1.3, which covers the case where 
the zeros of the characteristic function 
$$
f(z) = \E\,e^{izX}, \quad z  \in \C,
$$
of the subgaussian random variable $X$ are not necessarily real, 
but belong to the angle $|{\rm Arg}(z)| \leq \frac{\pi}{8}$.
Let us state it once more together with the stronger property (1.5).

\vskip5mm
{\bf Theorem 8.1.} {\sl Let $X$ be a subgaussian random variable 
with a symmetric distribution. If all zeros of $f(z)$ with
${\rm Re}(z) \geq 0$ lie in the angle $|{\rm Arg}(z)| \leq \frac{\pi}{8}$, 
then $X$ is strictly subgaussian. Moreover, if $X$ is not normal, 
then for any $t_0 > 0$, there exists $c = c(t_0)$,
$0 < c < \sigma^2 = \Var(X)$, such that
\be
\E\,e^{tX} \leq e^{c t^2/2}, \quad |t| \geq t_0.
\en
}

%\vskip2mm
In the proof of (8.1) we employ Proposition 2.5,
which asserts that (8.1) would follow from the property that
the function $t \rightarrow \log \E\,e^{\sqrt{t} X}$ is concave
on the positive half-axis $t \geq 0$ (in the symmetric case).
In this connection recall Proposition 5.1:
A random variable $\xi$ with characteristic function
$$
f_\xi(t) = e^{-t^2/2}\,(1 - \alpha t^2 + \beta t^4)
$$
is strictly subgaussian, if and only if $\beta \geq 0$ and
$\alpha \geq \sqrt{2\beta}$. In fact, the latter
description is also equivalent to the concavity of the function 
$$
t \rightarrow \log \E\,e^{\sqrt{t} \xi} =
-\frac{1}{2}\,t + \log (1 + \alpha t + \beta t^2), \quad t \geq 0.
$$
That is, we have:

\vskip5mm
{\bf Lemma 8.2.} {\sl Given $\alpha,\beta \geq 0$, the function
$Q(t) = \log (1 + \alpha t + \beta t^2)$ is concave in $t \geq 0$,
if and only if $\alpha \geq \sqrt{2\beta}$, and then the function
$R(t) = \alpha t - Q(t)$ is convex and non-decreasing.
}

\vskip5mm
Indeed, by the direct differentiation,
$$
R'(t) = 
\frac{(\alpha^2 - 2\beta) t + \alpha \beta t}{1 + \alpha t + \beta t^2},
\quad Q''(t) =  - \frac{(\alpha^2 - 2\beta) + 
2\alpha \beta t + 2 \beta^2 t^2}{(1 + \alpha t + \beta t^2)^2},
$$
from which the claim readily follows.

\vskip2mm
{\bf Proof of Theorem 8.1.} We may assume that $X$ is not normal.
By the symmetry assumption, with every zero $z = x+iy$,
we have more zeros $\pm x \pm iy$. So, one may arrange all zeros
in increasing order of their moduli and by coupling 
$\pm z_1, \pm \bar z_1,\dots$. Let us enumerate only
the zeros $z_n = x_n + iy_n$ lying in the quadrant 
$x_n \geq 0$, $y_n \leq 0$ and deal with $-z_n, \bar z_n, -\bar z_n$
as associated zeros. If $z_n$ is real, 
then we have only one associated zero $-z_n$. For simplicity
of notations, let us assume that all zeros are complex.

Since $X$ is subgaussian, the characteristic function $f(t)$ may be 
extended from the real line to the complex plane 
as an entire function satisfying
$$
|f(z)| \leq e^{b\, {\rm Im}(z)^2/2}, \quad z \in \C,
$$
for some constant $b \geq 0$. Therefore, $f$ is a ridge entire function 
of order $\rho \leq 2$ and of a finite type like in the strongly subgaussian 
case. Thus, Hadamard's theorem is applicable, with parameters 
$\rho \leq 2$ and $p \leq 2$. In this case, the representation (3.3) takes the form
$$
f(z) = e^{P(z)} \prod_{n \geq 1} 
G_p(z/z_n)\, G_p(z/\bar z_n)G_p(-z/z_n)\, G_p(-z/\bar z_n).
$$
Here, the genus of the canonical product satisfies $p \leq 2$, and 
$P(z)$ is a polynomial of degree at most 2 such that $P(0) = 0$. Thus,
putting in the sequel $w_n = \frac{1}{z_n} = a_n + b_n i$, we have
\be
f(z) = e^{i\beta z - \gamma z^2/2} \prod_{n \geq 1} \pi_{p,n}(z)
\en
for some $\beta,\gamma \in \C$, where
$$
\pi_{p,n}(z) = 
G_p(w_n z)\, G_p(-w_n z) G_p(\bar w_n z)\, G_p(-\bar w_n z).
$$
By the symmetry assumption, $f(-z) = f(z)$ for all $z \in \C$.
Since also $\pi_{p,n}(-z) = \pi_{p,n}(z)$, we conclude that 
$\beta = 0$. Put 
$$
\alpha_n = w_n^2 + \bar w_n^2 = 2(a_n^2 - b_n^2), \quad 
\beta_n =  |w_n|^4 = (a_n^2 + b_n^2)^2.
$$

There are three cases for the values of the genus, 
$p=0$, $p=1$, and $p=2$, for which
$$
G_0(u) = 1 - u, \quad
G_1(u) = (1 - u)\,e^u, \quad
G_2(u) = (1 - u)\,e^{u + \frac{u^2}{2}}.
$$
Since
$$
G_1(u) G_1(-u) = 1 - u^2 = G_0(u) G_0(-u) \quad {\rm and} \quad
G_2(-u) G_2(u) = (1 - u^2)\,e^{u^2},
$$
(8.2) is simplified to
\be
f(z) = e^{-\gamma z^2/2} \prod_{n \geq 1} Q_{p,n}(z),
\en
where
$$
Q_{0,n}(z) = Q_{1,n}(z) = (1 - w_n^2 z^2)(1 - \bar w_n^2 z^2) =
1 - \alpha_n z^2 + \beta_n z^4
$$
and
$$
Q_{2,n}(z) =
(1 - w_n^2 z^2)(1 - \bar w_n^2 z^2)\,e^{(w_n^2 + \bar w_n^2) z^2} =
(1 - \alpha_n z^2 + \beta_n z^4)\,e^{\alpha_n z^2}.
$$
These functions are real-valued for $z = t \in \R$, as well as
$f(t)$, by the symmetry assumption on the distribution of $X$. 
Hence, necessarily $\gamma \in \R$. Moreover, we have $\gamma \geq 0$, 
since otherwise $f(t)$ would not be bounded on the real axis $t \in \R$.

Since ${\rm Arg}(z_n) = -{\rm Arg}(w_n)$, we have 
${\rm Arg}(w_n) \leq \frac{\pi}{8}$, by the main angle hypothesis.
In particular, $a_n > b_n > 0$ so that $\alpha_n > 0$
(since $x_n > 0$, $y_n < 0$). As already noticed in the proof of Theorem 1.4, 
the angle hypothesis is equivalent to the relation $\alpha_n^2 \geq 2\beta_n$.

Applying (8.3) with $z = it$, $t \in \R$, we get that
\be
\E\,e^{tX} = e^{\gamma t^2/2} \prod_{n \geq 1} Q_{p,n}(it)
\en
with positive factors given by
$$
Q_{0,n}(it) = Q_{1,n}(it) = 1 + \alpha_n t^2 + \beta_n t^4, \qquad
Q_{2,n}(it) = (1 + \alpha_n t^2 + \beta_n t^4)\,e^{-\alpha_n t^2}.
$$
We have already observed in the proof of Proposition 5.1 that, 
by the angle hypothesis,
\be
1 + \alpha_n t^2 + \beta_n t^4 < e^{\alpha_n t^2}, \quad t > 0,
\en
so that $Q_{2,n}(it) < 1$.
Moreover, this inequality was strengthened by improving the constant
$\alpha_n$ in the exponent, provided that $t$ is bounded away
from zero. We will thus repeat some steps from the proof of Proposition 5.1.
However, formally, we need to consider the three cases separately
according to the three possible values of $p$.

{\bf Genus} $p=2$. By the very definition of the genus,
$$
\sum_{n \geq 1} |w_n|^3 = \sum_{n \geq 1} (a_n^2 + b_n^2)^{3/2}
 = \sum_{n \geq 1} \beta_n^{3/4}
< \infty.
$$
Since $Q_{2,n}(it) = 1 + O(\beta_n t^4)$ as $t \rightarrow 0$,
the product in (8.4) is absolutely convergent. Moreover, the right-hand 
side of (8.4) near zero is $1 + \gamma t^2 + O(t^3)$.
Hence, necessarily $\gamma = \sigma^2$, and (8.4) becomes
\be
\E\,e^{tX} = e^{\sigma^2 t^2/2} \prod_{n \geq 1} Q_{2,n}(it).
\en
Recalling the bound $Q_{2,n}(it) \leq 1$, we conclude that
\be
\E\,e^{tX} \leq e^{\sigma^2 t^2/2}, \quad t \in \R,
\en
which means that $X$ is strictly subgaussian.
For the second claim of the theorem, write
\be
\E\,e^{tX} = e^{V(t^2)},
\en
where
$$
V(s) = \frac{1}{2} \gamma s + \sum_{n \geq 1} \log Q_n(it) = 
\frac{1}{2} \sigma^2 s + \sum_{n \geq 1} \big[
\log(1 + \alpha_n s + \beta_n s^2) - \alpha_n s\big], \quad s \geq 0,
$$
and define 
$$
W(s)  = \frac{1}{2} \sigma^2 s - V(s) = \sum_{n \geq 1} R_n(s), \quad 
R_n(s) = \alpha_n s - \log(1 + \alpha_n s + \beta_n s^2).
$$
By Lemma 8.2, and using the assumption $\alpha_n^2 \geq 2\beta_n$,
all $R_n(s)>0$ for $s>0$, representing convex increasing functions. 
Hence, $W$ is a convex increasing function with $W(0) = 0$.
It remains to apply Proposition 2.5, and we obtain the property (8.1).

{\bf Genus} $p=1$.
By definition, the following sum converges
$$
\sum_{n \geq 1} |w_n|^2 = \sum_{n \geq 1} (a_n^2 + b_n^2)
 = \sum_{n \geq 1} \beta_n^{1/2} < \infty.
$$
Since 
$$
\alpha_n = 2(a_n^2 - b_n^2) \leq 2(a_n^2 + b_n^2) = 
2 \beta_n^{1/2},
$$
the product in (8.4) is convergent. Moreover, the right-hand 
side of (8.4) near zero is 
$$
1 + \frac{1}{2}\gamma t^2 + t^2 \sum_{n \geq 1} \alpha_n + O(t^3).
$$
Hence, necessarily 
$\frac{1}{2}\sigma^2 = \frac{1}{2}\gamma + \sum_{n \geq 1} \alpha_n$,
so that the characteristic function and the Laplace transform 
admit the same representation (8.6). As a result, since the
summation property defining the genus became stronger,
we immediately obtain (8.7) and its improvement (8.1)
using the previous step.

{\bf Genus} $p=0$.
By definition, the following sum converges
$$
\sum_{n \geq 1} |w_n| = \sum_{n \geq 1} (a_n^2 + b_n^2)^{1/2}
 = \sum_{n \geq 1} \beta_n^{1/4} < \infty.
$$
Since this assumption is stronger than the one of the previous
step, while $Q_{0,n} = Q_{1,n}$, we are reduced to the previous
step.
\qed

\vskip5mm
\section{{\bf Proof of Theorem 1.5}}
\setcounter{equation}{0}

\vskip2mm
\noindent
As in the proof of Theorem 8.1, let us enumerate the points 
$z_n = x_n + iy_n$ lying in the quadrant $x_n \geq 0$, $y_n \leq 0$ 
and deal with $-z_n, \bar z_n, -\bar z_n$ as associated zeros.
For simplicity of notations, we assume that all these numbers are complex. Put
$w_n = \frac{1}{z_n} = a_n + b_n i$ and define
\bee
f_n(z)
 & = &
e^{- \gamma_n z^2/2}\,
(1 - w_n z)(1+w_n z)(1 - \bar w_n z)(1 + \bar w_n z) \\
 & = &
e^{-\gamma_n z^2/2}\,(1 - \alpha_n z^2 + \beta_n z^4), \quad z \in \C,
\ene
for a given sequence $\gamma_n > 0$ (to be precised later on) with
$\alpha_n = 2(a_n^2 - b_n^2)$ and $\beta_n =  (a_n^2 + b_n^2)^2$ as before.
By the assumption, $a_n,b_n>0$. Moreover, the angle assumption
$|{\rm Arg}(z_n)| = {\rm Arg}(w_n) \leq \frac{\pi}{8}$ is
equivalent to $\alpha_n^2 \geq 2\beta_n$, which may also be
written as 
\be
b_n \leq \frac{1}{\sqrt{2} + 1}\,a_n.
\en

Now, if $\gamma_n$ is sufficiently large, $f_n(t)$, $t \in \R$, 
will be the characteristic function of a strictly subgaussian
distribution. A full description of the minimal possible value
of $\gamma_n$ is provided in Proposition 7.1. More precisely,
consider the function
\bee
g_n(t) \, = \, f_n\Big(\frac{t}{\sqrt{\gamma_n}}\Big)
 & = &
e^{- t^2/2}\,
(1 - w_n' t)(1+w_n' t)(1 - \bar w_n' t)(1 + \bar w_n' t) \\
 & = &
e^{-t^2/2}\,(1 - \alpha_n' t^2 + \beta_n' t^4)
\ene
with
$$
w_n' = a_n' + b_n' i, \quad 
a_n' = \frac{a_n}{\sqrt{\gamma_n}}, \quad 
b_n' = \frac{b_n}{\sqrt{\gamma_n}}, \quad
\alpha_n' = \frac{2(a_n^2 - b_n^2)}{\gamma_n}, \quad 
\beta_n' = \frac{(a_n^2 + b_n^2)^2}{\gamma_n^2}.
$$
As we know, $g_n(t)$ represents the characteristic function of 
a strictly subgaussian random variable $X_n'$, as long as
$$
b_n' \leq \frac{1}{\sqrt{2} + 1}\,a_n', \quad a_n' \leq a_0,
$$
where the universal constant $a_0$ was explicitly identified
in (7.6), $a_0 \sim 0.7391$. Here, the first condition is satisfied
in view of (9.1), while the second one is equivalent to
\be
\gamma_n \geq \frac{a_n^2}{a_0^2}.
\en
Moreover, $X_n'$ has variance
$$
\Var(X_n') = - g_n''(0) = 2\alpha_n' + 1 = 
\frac{4(a_n^2 - b_n^2)}{\gamma_n} + 1.
$$
Thus, subject to (9.2), $f_n(t)$ will be the characteristic 
function of the strictly subgaussian random variable 
$X_n = \sqrt{\gamma_n}\, X_n'$, whose variance is given by
\be
\Var(X_n) =  4(a_n^2 - b_n^2) + \gamma_n.
\en

Now, assuming that $\Lambda \geq 4 + \frac{1}{a_0^3} \sim 5.83$,
let us choose
$$
\gamma_n = (\Lambda - 4) a_n^2 + (\Lambda + 4) b_n^2,
$$
so that the expression in (9.3) would be equal to
$\Lambda (a_n^2 + b_n^2) = \Lambda |w_n|^2$. Then
the condition (9.2) is satisfied, and also $\sum_n \gamma_n < \infty$.
As a result, the series $\sum_n X_n$ is convergent with 
probability one, and the sum of the series, call it $X$, 
represents a strictly subgaussian random variable with 
characteristic function
$$
f(z) = \prod_n f_n(z)
$$
(cf. Proposition 2.2). By the construction, all
$f_n(z)$ have exactly prescribed zeros, and
$$
\Var(X) = \sum_n \Var(X_n) = \Lambda \sum_n |w_n|^2.
$$
\qed

\vskip5mm
\section{{\bf Laplace Transforms with Periodic Components}}
\setcounter{equation}{0}

\vskip2mm
\noindent
We now turn to a second class of Laplace transforms -- the ones that
contain periodic components. Recall that a random variable $X$  belongs 
to the class $\mathfrak F_h$, $h>0$, if it has a density $p(x)$ such that
the function
$$
q(x) = \frac{p(x)}{\varphi(x)},
\quad x \in \R,
$$
is periodic with period $h$. This section is devoted to basic properties
of this class (some of them will be used in the proof of Theorem 1.6).

\vskip5mm
{\bf Proposition 10.1.} {\sl If $X$ belongs to the class 
$\mathfrak F_h$, then for all integers $m$,
\be
\E\,e^{mh X} = e^{(mh)^2/2}.
\en
In particular, the random variable $X$ is subgaussian. 
}

\vskip5mm
{\bf Proof.}
By the periodicity of $q$, the random variable $X + mh$ has density
\bee
p(x - mh) 
 & = &
q(x - mh) \varphi(x - mh) \\
 & = &
q(x)\varphi(x)\,e^{mh x - (mh)^2/2} \, = \, 
p(x)\,e^{mh x - (mh)^2/2}.
\ene
It remains to integrate this equality over $x$, which leads to (10.1).

Next, starting from (10.1), it is easy to see that $\E\,e^{c X^2} < \infty$
for some $c>0$.
\qed

\vskip5mm
As a consequence, the Laplace transform $L(t) = \E\,e^{tX}$, $t \in \R$,
is finite and may be extended to the complex plane as an entire function.
This property may be refined.

\vskip5mm
{\bf Proposition 10.2.} {\sl If $X$ belongs to $\mathfrak F_h$, then
its Laplace transforms is an entire function of order 2.
Moreover, if $\E X = 0$, it satisfies
\be
|L(z)| \leq e^{(|t|+h)^2/2}, \quad z = t+iy \in \C.
\en
}

{\bf Proof.} We may assume that $\E X = 0$. In this case,
by Jensen's inequality, $L(t) \geq 1$ for all $t \in \R$, so that
$t=0$ is the point of miminum of $L$ on the real line.
Since $L(t)$ is convex (and moreover, $\log L(t)$ is convex),
$L(t)$ is decreasing for $t<0$ and is increasing for $t>0$.

Given $t \geq 0$, take an integer number $m \geq 1$ such that
$(m-1) h \leq t < mh$. Then, by (10.1), and using the monotonicity of $L$,
we get
\be
L(t) \leq L(mh) = e^{(mh)^2/2} \leq  e^{(t+h)^2/2}.
\en
By a similar argument, $L(-t) \leq e^{(t+h)^2/2}$.
Thus, we obtain (10.2) for real values of $z$ (when $y=0$).
In the general case, it remains just to note that
$|L(z)| \leq L(t)$, and we obtain (10.2). This bound shows that
$L(z)$ is an entire function of order at most 2.

On the other hand, (10.1) shows that $L(z)$ is an entire function of order at least 2.
\qed

\vskip5mm
{\bf Proposition 10.3.} {\sl If $X$ belongs to $\mathfrak F_h$, then
the function
\be
\Psi(t) = L(t)\,e^{-t^2/2}, \quad t \in \R,
\en
is periodic with period $h$. It can be extended to the complex plane
as an entire function. Moreover, if $\E X = 0$, it
satisfies
\be
|\Psi(z)| \leq C_{h,y}\,e^{h|t|}, \quad z = t+iy \in \C,
\en
with $C_{h,y} = e^{(h^2 + y^2)/2}$.
}

\vskip2mm
The inequality (10.5) shows that $\Psi$ is an entire function
of order at most 2. 

By analyticity and periodicity on the real line,
\be
\Psi(z+h) = \Psi(z) \quad {\rm for \ all} \ z \in \C.
\en

\vskip2mm
{\bf Proof.} By periodicity of $q$, changing the variable $x = y+h$,
we have
\bee
L(t+h) 
 & = &
\int_{-\infty}^\infty e^{(t+h)\,x}\,q(x)\,\varphi(x)\,dx 
 \, = \,
\int_{-\infty}^\infty e^{(t+h)\,(y+h)}\,q(y+h)\,\varphi(y+h)\,dy \\
 & = &
\int_{-\infty}^\infty e^{(t+h)\,(y+h)}\,q(y)\,\varphi(y)\,e^{-yh - h^2/2}\,dy
 \, = \,
L(t)\,e^{th + h^2/2}.
\ene
Hence
$$
L(t+h) \,e^{-(t+h)^2/2} = L(t)\,e^{-t^2/2},
$$ 
which was the first claim. Since $L(z)$ is an entire function,
$\Psi(z)$ is entire as well.

Next, assuming that $\E X = 0$, one may apply (10.2) which gives 
$\Psi(t) \leq C_h\,e^{h |t|}$ with  $C_h = e^{h^2/2}$. Thus, we 
obtain (10.5) for real values of $z$. In the general case, for simplicity 
let $t = {\rm Re}(z) \geq 0$. By the previous step, 
$$
|L(z)| \leq L(t) \leq C_h\,e^{ht}.
$$
Hence
$$
|\Psi(z)| \leq L(t) \,e^{{\rm Re}(z^2)/2} \leq 
e^{(t+h)^2/2}\,e^{(y^2 - t^2)/2} = C_{h,y}\,e^{ht}.
$$
\qed

\vskip2mm
Let us also examine the periodicity property for convolutions.
Here, the basic observation concerns the normalized sums
$Z_n = \frac{X_1 + \dots + X_n}{\sqrt{n}}$,
where $X_k$'s are independent copies of the random variable $X$.

\vskip5mm
{\bf Proposition 10.4.} {\sl If $X$ belongs to $\mathfrak F_h$, then
$Z_n$ belongs to $\mathfrak F_{h\sqrt{n}}$.
}

\vskip5mm
{\bf Proof.} Let $n=2$ for simplicity
of notations. Let $p$ be the density of $X$ such that $q = p/\varphi$
is $h$-periodic. Since $X_1+X_2$ has density
$$
\int_{-\infty}^\infty p(x - y)\, p(y)\,dy =
\int_{-\infty}^\infty p\Big(\frac{x}{2} + z\Big)\, p\Big(\frac{x}{2} - z\Big)\,dz,
$$
the density of $Z_2 = \frac{X_1 + X_2}{\sqrt{2}}$ may be written as
\bee
p_2(x)
 & = &
\sqrt{2} \int_{-\infty}^\infty q\Big(\frac{x}{\sqrt{2}} + z\Big)\, 
q\Big(\frac{x}{\sqrt{2}} - z\Big)\,\varphi\Big(\frac{x}{\sqrt{2}} + z\Big)\, 
\varphi\Big(\frac{x}{\sqrt{2}} - z\Big)\,dz \\
& = &
\frac{1}{\sqrt{\pi}}\,\varphi(x)
\int_{-\infty}^\infty q\Big(\frac{x}{\sqrt{2}} + z\Big)\, 
q\Big(\frac{x}{\sqrt{2}} - z\Big)\,e^{-z^2}\,dz.
\ene
Thus, the correspondong $q$-function for $Z_2$ is given by
$$
q_2(x) = \frac{p_2(x)}{\varphi(x)} = \frac{1}{\sqrt{\pi}}
\int_{-\infty}^\infty q\Big(\frac{x}{\sqrt{2}} + z\Big)\, 
q\Big(\frac{x}{\sqrt{2}} - z\Big)\,e^{-z^2}\,dz.
$$
As $q$ is $h$-periodic, the last integrand is periodic with respect to
the variable $x$, with period $h\sqrt{2}$. Consequently, 
$q_2(x+h\sqrt{2}) = q_2(x)$ for all $x \in \R$.
\qed

\vskip5mm
\section{{\bf Proof of Theorem 1.6}}
\setcounter{equation}{0}

\vskip2mm
\noindent
In view of the previous observations, we only need to consider
the necessity part in the statement of Theorem 1.6 and prove
the periodicity of the density $q$.

Since $X$ is subgaussian, its Laplace transform
is an entire function of order at most 2. Hence $\Psi(z)$ is
also entire and satisfies (10.6). Thus, the characteristic function 
of $X$ is an entire function representable in the complex plane as
$f(z) = \Psi(iz)\,e^{-z^2/2}$. Hence, by (10.6),
$$
f(t+ih)\,e^{(t+ih)^2/2} = f(t)\,e^{t^2/2}
$$
for all $t \in \R$, that is,
\be
f(t+ih) = f(t)\,e^{-ith + h^2/2}.
\en

By the integrability assumption, the random variable $X$ 
has a continuous density $p(x)$ given by the Fourier inversion formula
$$
p(x) = 
\frac{1}{2\pi} \int_{-\infty}^\infty e^{-itx} f(t)\,dt, \quad x \in \R.
$$
This yields
$$
q(x) = \frac{p(x)}{\varphi(x)} = 
\frac{1}{\sqrt{2\pi}}\,e^{x^2/2} \int_{-\infty}^\infty e^{-itx} f(t)\,dt
$$
and
$$
q(x+h) = \frac{1}{\sqrt{2\pi}}\,e^{x^2/2}\, e^{xh + h^2/2} 
\int_{-\infty}^\infty e^{-itx - ith} f(t)\,dt.
$$
Hence, we need to show that
\be
\int_{-\infty}^\infty e^{-itx} f(t)\,dt = e^{xh + h^2/2} 
\int_{-\infty}^\infty e^{-itx - ith} f(t)\,dt.
\en

Using contour integration, one may rewrite the first integral
in a different way. Given $T>0$, consider the rectangle contour with sides
\bee
C_1 & = &
[-T,T], \qquad \qquad \quad C_2 = [T,T+ih], \\ 
C_3 & = &
[T+ih,-T+ih], \quad C_4 = [-T+ih,-T],
\ene
so that to apply Cauchy's theorem and write down
$$
\int_{C_1} e^{-izx} f(z)\,dz + \int_{C_2} e^{-izx} f(z)\,dz +
\int_{C_3} e^{-izx} f(z)\,dz + \int_{C_4} e^{-izx} f(z)\,dz = 0.
$$
For points $z = t+iy$ on the contour, we have
$|e^{-izx}| = e^{xy} \leq e^{|x|h}$. In addition, $f(z) \rightarrow 0$
as $|t| \rightarrow \infty$ uniformly over all $y$ such that
$|y| \leq h$. This follows from the fact that the functions
$t \rightarrow f(t+iy)$ represent the Fourier transform of
the functions $p_y(x) = e^{-xy} p(x)$. Indeed, by the subgaussian
assumption, the family $\{p_y: |y| \leq h\}$ is pre-compact in
$L^1(\R^n)$, so that the Riemeann-Lebesgue lemma is applicable
to the whole family. As a consequence,
\bee
\int_{-\infty}^\infty e^{-itx} f(t)\,dt
 & = &
\lim_{T \rightarrow \infty} \int_{C_1} e^{-izx} f(z)\,dz \\
 & & \hskip-10mm = \ -
\lim_{T \rightarrow \infty} \int_{C_3} e^{-izx} f(z)\,dz \, = \,
\int_{-\infty}^\infty e^{-i(t + ih)x} f(t+ih)\,dt.
\ene
where the last integral is convergent due to (11.1).
Moreover, by (11.1), the last integrand is equal to
$$
e^{-i(t + ih)x}\, \,e^{-ith + h^2/2}\,f(t),
$$
which coincides with the integrand on the right-hand side
of (11.2) multiplied by the indicated factor. This proves (11.2).
\qed

\vskip5mm
{\bf Remark 11.1.} Since $f(t) = L(it) = \Psi(it)\,e^{-t^2/2}$,
the integrability assumption in Theorem 1.6 is fulfilled,
as long as $\Psi(z)$ has order smaller than 2.

\vskip5mm
\section{{\bf Examples Involving Triginometric Series}}
\setcounter{equation}{0}

\vskip2mm
\noindent
Theorem 1.6 is applicable to a variety of interesting examples including
the underlying distributions whose Laplace transform has the form 
$$
L(t) = \Psi(t)\, e^{t^2/2}, \quad t \in \R,
$$ 
where $\Psi$ is a $2\pi$-periodic functions of the form
\be
\Psi(t) = 1 - c P(t), \quad P(t) = a_0 + \sum_{k=1}^\infty (a_k \cos(kt) + b_k \sin(kt)).
\en
Here $a_k,b_k$ are real coefficients which are supposed to satisfy
\be
\sum_{k=1}^\infty e^{k^2/2}\,(|a_k| + |b_k|) < \infty,
\en
and $c \in \R$ is a non-zero parameter.

\vskip5mm
{\bf Proposition 12.1.} {\sl If $P(0) = P'(0) = P''(0) = 0$ and $|c|$ is small 
enough, then $L(t)$ represents the Laplace transform of a subgaussian 
random variable $X$ with $\E X = 0$, $\E X^2 = 1$, and with density 
$p = q \varphi$, where $q$ is a bounded, $2\pi$-periodic function. This 
random variable is strictly subgaussian, if $P(t) \geq 0$ for all $t \in \R$ 
and if $c > 0$ is small enough.
}

\vskip5mm
{\bf Proof.} The functions of the form
$u_\lambda(x) = \cos(\lambda x)\,\varphi(x)$ and 
$v_\lambda(x) = \sin(\lambda x)\,\varphi(x)$
have respectively the Laplace transforms
\bee
\int_{-\infty}^\infty e^{tx} u_\lambda(x)\,dx
 & = &
e^{-\lambda^2/2}\,\cos(\lambda t)\,e^{t^2/2}, \\
\int_{-\infty}^\infty e^{tx} v_\lambda(x)\,dx
 & = &
e^{-\lambda^2/2}\,\sin(\lambda t)\,e^{t^2/2}.
\ene
Define
\be
q(x) = \varphi(x) - c \varphi(x)\Big(a_0 +
\sum_{k=1}^\infty e^{k^2/2}\,(a_k \cos(kx) + b_k \sin(kx))\Big).
\en
In this case, the Laplace transform of the function 
$p(x) = q(x) \varphi(x)$ is exactly
$$
\int_{-\infty}^\infty e^{tx} p(x)\,dx =
(1 - cP(t))\,e^{-t^2/2}.
$$

The requirement $P(0) = 0$ guarantees that $\int_{-\infty}^\infty p(x)\,dx = 1$. 
Moreover, according to (12.3), the condition on the parameter $c$ which ensures 
that the function $p$ is indeed a probability density may be stated as
$$
|a_0| + \sum_{k=1}^\infty e^{k^2/2}\,(|a_k| + |b_k|) \leq \frac{1}{|c|}.
$$
This is fulfilled due to (12.1) when $|c|$ is small enough. Finally, the properties 
$\E X = 0$, $\E X^2 = 1$ are equivalent to $P'(0) = P''(0) = 0$.
\qed

\vskip2mm
Note that in terms of the coeficients in the series (12.1), the condition
$P(0) = P'(0) = P''(0) = 0$ has the form
$$
a_0 + \sum_{k=1}^\infty a_k = \sum_{k=1}^\infty k b_k = \sum_{k=1}^\infty k^2 a_k = 0.
$$

It should also be mentioned that, when $P$ is a trigonometric polynomial of degree $N$, 
the function $q(x)$ in (12.3) will be a trigonometric polynomial of degree $N$ as well.

\vskip5mm
{\bf Example 12.2.} As a particular case, one may consider the transforms
\be
L(t) = (1 - c \sin^m(t))\, e^{t^2/2}
\en
with an arbitrary integer $m \geq 3$, where $c$ is small enough.
Then $\E X = 0$, $\E X^2 = 1$, and the cumulants of $X$ satisfy
$$
\gamma_k(X) = 0, \quad 3 \leq k \leq m-1.
$$

Moreover, if $m \geq 4$ is even and $c>0$, the random variable 
$X$ with the Laplace transform (12.3) is strictly subgaussian.
In the case $m=4$, (12.4) corresonds to the $\pi$-periodic 
polynomial $P(t) = \sin^4 t = \frac{1}{8}\,(3 - 4 \cos(2t) + \cos(4t))$.

\vskip5mm
\section{{\bf Examples Involving Poisson Formula and Theta Functions}}
\setcounter{equation}{0}

\vskip2mm
\noindent
Often, the periodic functions $\Psi(t)$ in (12.1) appear naturally by means
of the Poisson formula, rather than as a trigonometric series.
Let $w(t) \geq 0$ be an integrable, even, absolutely continuous function 
on the real line with Fourier transform
$$
\hat w(x) = \int_{-\infty}^\infty e^{itx} w(t)\,dt, \quad x \in \R.
$$
As a natrural generalization of Example 12.2 with $m=4$, we have
the following corollary from Proposition 12.1 assuming that
\be
\sum_{k=1}^\infty e^{(k+4)^2/2}\,|\hat w(k)| < \infty.
\en

\vskip2mm
{\bf Corollary 13.1.} {\sl For all $c>0$ small enough,
$$
L(t) = \Psi(t)\, e^{t^2/2}, \quad \Psi(t) = 1 - c \,(\sin t)^4  \sum_{m \in \Z} w(t + 2\pi m),
$$
represents the Laplace transform of a strictly subgaussian random variable $X$ with $\E X = 0$, 
$\E X^2 = 1$, which has density $p(x) = q(x) \varphi(x)$, where $q(x)$ is a $2\pi$-periodic function.
}

\vskip5mm
{\bf Proof.} The function $Q(t) = \sum_{m \in \Z} w(t + 2\pi m)$
is well-defined (since the series is absolutely convergent), $2\pi$-periodic, 
and admits a Fourier series expansion
$$
Q(t) = \frac{1}{2\pi} \sum_{k \in \Z} \hat w(k)\, e^{-ikt}.
$$
This is a well-known Poisson formula, in which the series is understood
as a limit of symmetric partial sums, cf. e.g. \cite{Z}, p. 68. Under (13.1), this series
is absolutely convergent and defines a smooth function. By the symmetry 
$\hat w(-k) = \hat w(k)$, $k \in \Z$, this formula takes the form
$$
Q(t) = \frac{1}{2\pi}\,\Big[\hat w(0) + 2\sum_{k = 1}^\infty \hat w(k) \cos(kt)\Big].
$$

Using $\sin^4 t  = \frac{1}{8}\,(3 - 4 \cos(2t) + \cos(4t))$, we have
\bee
16 \pi\, (\sin t)^4 \, Q(t)
 & = &
3\hat w(0)  - 4 \hat w(0) \cos(2t) + \hat w(0) \cos(4t) \\
 & & \hskip-25mm + \
6\sum_{k = 1}^\infty \hat w(k) \cos(kt) - 
8 \sum_{k = 1}^\infty \hat w(k) \cos(kt) \cos(2t) +  
2 \sum_{k = 1}^\infty \hat w(k) \cos(kt) \cos(4t).
\ene
Applying the identity
$\cos a \cos b = \frac{1}{2}\,\cos(a+b) + \frac{1}{2}\,\cos(a-b)$, one may
rewrite the last line~as
$$
6\sum_{k = 1}^\infty \hat w(k) \cos(kt) - 
4 \sum_{k = 1}^\infty \hat w(k) (\cos((k+2) t) + \cos((k-2)\,t))
$$
$$
+ \sum_{k = 1}^\infty \hat w(k) (\cos((k+4) t) + \cos((k-4)\,t)).
$$
Hence for $k \geq 5$ the coefficients $a_k$ in the Fourier series for
$P(t) = (\sin t)^4\, Q(t)$ are given by
$$
16 \pi\,a_k = 6 \hat w(k) - 4 \hat w(k+2) - 4 \hat w(k-2) + \hat w(k+4) + \hat w(k-4).
$$
Hence, the condition (12.2) is fulfilled under (13.1), and one may apply Proposition 12.1.
\qed

\vskip5mm
{\bf Example 13.2.} One may further apply Corollary 13.1 to the theta functions
$Q(t)$ corresponding to
$$
w(t) = \frac{1}{\sigma \sqrt{2\pi}}\,e^{-t^2/2\sigma^2}, \quad
\hat w(x) = e^{-\sigma^2 x^2/2}
$$
with an arbitrary parameter $\sigma > 1$.

\vskip5mm
\section{{\bf Central Limit Theorem for R\'enyi Distances}}
\setcounter{equation}{0}

\vskip2mm
\noindent
Finally, let us describe the role of subgaussian distributions in the 
central limit theorem with respect to the R\'enyi divergences $D_\alpha$
defined in (1.3). Consider the normalized sums
$$
Z_n = \frac{X_1 + \dots + X_n}{\sqrt{n}},
$$
where $X_k$' are independent copies of a random variable $X$ with
mean zero and variance one. Assuming that $Z_n$ have densities $p_n$
for some or equivalently for all sufficiently large $n$,
the following characterization was obtained in \cite{B-C-G1},
which we state in dimension one.

\vskip5mm
{\bf Theorem 14.1.} {\sl Fix $1 < \alpha  < \infty$. For the convergence 
\be
D_\alpha(p_n||\varphi) \rightarrow 0 \ {\rm as} \ n \rightarrow \infty, 
\en
it is necessary 
and sufficient that $D_\alpha(p_n||\varphi) <\infty$ for some $n = n_0$, and
\be
\E\,e^{tX} < e^{\beta t^2/2} \quad {\sl for \ all} \ t \neq 0,
\en
where $\beta = \frac{\alpha}{\alpha-1}$ is the conjugate index.
}

\vskip5mm
Thus, for the CLT as in (14.1), the random variable $X$ has to be subgaussian.
In order to obtain this convergence for all $\alpha$ simultaneously, the
condition (14.2) on the Laplace transform should be fulfilled for all $\beta > 1$. 
But this is equivalent to saying that $X$ is strictly subgaussian, thus proving
Theorem 1.1.

In this connection, it is natural to raise the question of whether or not (14.1)
may hold for the critical index $\alpha = \infty$, which corresponds to the
strongest distance in this hierarchy. Note that in the limit case it is defined to be
$$
D_\infty(p_n||\varphi) = \lim_{\alpha \rightarrow \infty} D_\alpha(p_n||\varphi) =
\log \, {\rm ess\,sup}_{x \in \R} \ \frac{p_n(x)}{\varphi(x)}.
$$
As an equivalent quantity, one may also consider the limit Tsallis distance
$$
T_\infty(p_n||\varphi) =
{\rm ess\,sup}_{x \in \R} \ \frac{p_n(x) - \varphi(x)}{\varphi(x)}.
$$
Suppose it is finite for some $n = n_0$. The following two theorems
can be obtained using recent results on the sharpened Richter-type local 
limit theorem, cf. \cite{B-C-G2}.

\vskip5mm
{\bf Theorem 14.2.} {\sl Suppose that, for every $t_0>0$,
\be
\E\,e^{tX} \leq \delta e^{t^2/2} \quad {\sl for \ all} \ |t| \geq t_0
\en
with some $\delta =\delta(t_0) \in (0,1)$. Then
\be
T_\infty(p_n||\varphi) = O\Big(\frac{(\log n)^3}{n}\Big) \quad {\sl as} \
n \rightarrow \infty.
\en
}

Note that (14.3) is a weakened form of the separation property (1.5), which
in turn is a sharpening of strict subgaussianity.
In particular, this rate for the convergence in $D_\infty$ holds true for all
distributions from the class $\mathfrak L$ whose densities $p(x)$ are dominated
by $\varphi(x)$.

A similar assertion holds true in the period case.

\vskip5mm
{\bf Theorem 14.3.} {\sl Suppose that $X$ is strictly subgaussian, with an
$h$-periodic function $\Psi(t) = L(t)\,e^{-t^2/2}$, $h>0$. If $\Psi(t) < 1$ 
in the interval $0 < t < h$, then $(14.4)$ is true as well.
}

\end{document}